\documentclass[12pt]{article} 
 \usepackage{ams}
 \usepackage{dsfont}
 \usepackage{amsmath}
 \usepackage{amsthm}
 \usepackage{txfonts,bm}
 \usepackage{bold-extra}
 \usepackage[sc,bf]{caption}
 \usepackage{rotating}
 \usepackage{fancyhdr}
 \usepackage{fancyheadings}
 \RequirePackage[OT1]{fontenc}
 \RequirePackage{amsthm,amsmath}
 \RequirePackage[numbers]{natbib}
 \usepackage[colorlinks,citecolor=blue,urlcolor=blue,linkcolor=blue]{hyperref}

 \usepackage{titlesec}
 \usepackage{latexsym}
   \font\twelvebm                       = cmmib10 at 12truept
   \font\tenbm                          = cmmib10 at 10truept
   \font\sevenbm                        = cmmib10 at 7truept
 = \tenbm \scriptscriptfont9              = \sevenbm
 
   at 10truept  at 10truept  at
 10truept  at 10truept
  at 10truept

 \mathchardef \BGamma            = "0900 \mathchardef \BDelta
 = "0901 \mathchardef \BTheta            = "0902 \mathchardef
 \BLambda           = "0903 \mathchardef \BXi               = "0904
 \mathchardef \BPi               = "0905 \mathchardef \BSigma
 = "0906 \mathchardef \BUpsilon          = "0907 \mathchardef \BPhi
 = "0908 \mathchardef \BPsi              = "0909 \mathchardef
 \BOmega            = "090A \mathchardef \Balpha            = "090B
 \mathchardef \Bbeta             = "090C \mathchardef \Bgamma
 = "090D \mathchardef \Bdelta            = "090E \mathchardef
 \Bepsilon          = "090F \mathchardef \Bzeta             = "0910
 \mathchardef \Beta              = "0911 \mathchardef \Btheta
 = "0912 \mathchardef \Biota             = "0913 \mathchardef
 \Bkappa            = "0914 \mathchardef \Blambda           = "0915
 \mathchardef \Bmu               = "0916 \mathchardef \Bnu
 = "0917 \mathchardef \Bxi               = "0918 \mathchardef \Bpi
 = "0919 \mathchardef \Brho              = "091A \mathchardef
 \Bsigma            = "091B \mathchardef \Btau              = "091C
 \mathchardef \Bupsilon          = "091D \mathchardef \Bphi
 = "091E \mathchardef \Bchi              = "091F \mathchardef \Bpsi
 = "0920 \mathchardef \Bomega            = "0921 \mathchardef
 \Bvarepsilon       = "0922 \mathchardef \Bvartheta         = "0923
 \mathchardef \Bvarpi            = "0924 \mathchardef \Bvarrho
 = "0925 \mathchardef \Bvarsigma         = "0926 \mathchardef
 \Bvarphi           = "0927
 \mathchardef \bA        = "0941 \mathchardef \bB        = "0942
 \mathchardef \bC        = "0943 \mathchardef \bD        = "0944
 \mathchardef \bE        = "0945 \mathchardef \bF        = "0946
 \mathchardef \bG        = "0947 \mathchardef \bH        = "0948
 \mathchardef \bI        = "0949 \mathchardef \bJ        = "094A
 \mathchardef \bK        = "094B \mathchardef \bL        = "094C
 \mathchardef \bM        = "094D \mathchardef \bN        = "094E
 \mathchardef \bO        = "094F \mathchardef \bP        = "0950
 \mathchardef \bQ        = "0951 \mathchardef \bR        = "0952
 \mathchardef \bS        = "0953 \mathchardef \bT        = "0954
 \mathchardef \bU        = "0955 \mathchardef \bV        = "0956
 \mathchardef \bW        = "0957 \mathchardef \bX        = "0958
 \mathchardef \bY        = "0959 \mathchardef \bZ        = "095A
 \mathchardef \ba        = "0961 \mathchardef \bb        = "0962
 \mathchardef \bc        = "0963 \mathchardef \bd        = "0964
 \mathchardef \bee       = "0965 
 \mathchardef \bff       = "0966 \mathchardef \bg        = "0967
 \mathchardef \bh        = "0968
 \mathchardef \bj        = "096A \mathchardef \bk        = "096B
 \mathchardef \bl        = "096C \mathchardef \bm        = "096D
 \mathchardef \bn        = "096E \mathchardef \bo        = "096F
 \mathchardef \bp        = "0970 \mathchardef \bq        = "0971
 \mathchardef \br        = "0972 \mathchardef \bs        = "0973
 \mathchardef \bt        = "0974 \mathchardef \bu        = "0975
 \mathchardef \bv        = "0976 \mathchardef \bw        = "0977
 \mathchardef \bx        = "0978 \mathchardef \by        = "0979
 \mathchardef \bz        = "097A

 \font\tencb            = cmssbx10 scaled \magstep4 \font\eigcb
 = cmssbx10 scaled \magstep2 \textfont8             = \tencb
 \mathchardef\bAs       = "1841
 \def\Asem#1#2{\mathop{\vrule height10.5pt depth5.5pt width0pt\bAs}_{#1}^{#2}}
 \def\asem#1#2{
          \ifmmode
         \ifinner
            \raise0.9pt\hbox{$\scriptstyle\bAs$}_{#1}^{#2}
         \else
            \Asem{#1}{#2}
         \fi
          \fi
          }

 \newtheorem{theo}{\small\bf Theorem}[section]
 \newtheorem{lem}{\small\bf Lemma}[section]
 
 \newtheorem{rem}{\small\bf Remark}[section]
 
 \newtheorem{exam}{\small\bf Example}[section]
 
  \newtheorem{defi}{\small\bf Definition}[section]

 \renewcommand{\Pr}{\mbox{\rm  \hspace*{.2ex}I\hspace{-.5ex}P\hspace*{.2ex}}}

 \newcommand{\be}{\begin{equation}}
 \newcommand{\ee}{\end{equation}}

 \newcommand{\E}{\mbox{\rm \hspace*{.2ex}I\hspace{-.5ex}E\hspace*{.2ex}}}

 \newcommand{\Var}{\mbox{\rm \hspace*{.2ex}Var\hspace*{.2ex}}}
 
 \newcommand{\Cov}{\mbox{\rm \hspace*{.2ex}Cov\hspace*{.2ex}}}

 \newcommand{\sign}{\mbox{sign}}

 \newcommand{\bbb}[1]{\mbox{\boldmath $ #1 $}}

 \newenvironment{pr}[1]{{\small\bf {#1}:}}{}

 \newcommand{\RR}{\mathbb{R}}

 \numberwithin{equation}{section}

 \title{\Large\bf A discrete analogue of Terrell's characterization
 of rectangular distributions}

 \author{\large
 Nickos
 Papadatos
 {\footnote{
 {e-mail:}\
 \textcolor[rgb]{0.98,0.00,0.00}{npapadat@math.uoa.gr}, {url:}\
 \textcolor[rgb]{0.98,0.00,0.00}{users.uoa.gr/$\sim$npapadat/}}}}
 
 \vspace{-1em}

 \date{\small\it
 \begin{tabular}{r@{\hspace{0ex}}l}
 & National and Kapodistrian University of Athens, Department of Mathematics,
 \\
 [-.3ex]
 &
 Section of Statistics and Operations
 Research,
 Panepistemiopolis, 157 84 Athens, Greece.
 \end{tabular}
 }
 \vspace{-1em}

\begin{document}

 \maketitle
 \vspace*{-2em}


 \thispagestyle{empty}

 \begin{abstract}
 \noindent
 George R.\ Terrell (1983, {\it Ann.\ Probab.}, vol.\ 11(3), pp.\ 823--826)
 showed that the Pearson coefficient of correlation of an ordered pair from a
 random sample of size two is at most one-half, and the equality is attained
 only for rectangular (uniform over some interval) distributions.

 In the present note it is proved that the same is true for the discrete case,
 in the sense that the correlation coefficient attains its maximal value
 only for discrete rectangular (uniform over some finite lattice) distributions.
 \end{abstract}
 {\footnotesize {\it MSC}:  Primary 60E15; 62E10; Secondary 62G30.
 }
 \vspace*{-1em}

 \section{Introduction and main result}
 \label{sec.1}
 For independent, identically distributed, random variables $X_1$, $X_2$,
 from a probability distribution function $F$, the corresponding order statistics
 will be denoted by $X_{1:2}\leq X_{2:2}$, that is, $X_{1:2}=\min\{X_1,X_2\}$,
 $X_{2:2}=\max\{X_1,X_2\}$, and the Pearson coefficient of correlation
 by
 \[
 \rho_{12}:=\frac{\Cov(X_{1:2},X_{2:2})}{\sqrt{\Var X_{1:2}}\sqrt{\Var X_{2:2}}}.
 \]
 For a random pair $(X,Y)$, the correlation coefficient,
 $\rho(X,Y)$, is well defined (and belongs to the interval $[-1,1]$)
 if and only if both $X$, $Y$ are non-degenerate with finite second moment.
 Thus, for $\rho_{12}$ to be well defined it is necessary and sufficient that $F$ is non-degenerate and posses finite second moment (so that
 $0<\Var X_{j:2}<\infty$, $j=1,2$).

 It is known for a long time
 that $\rho_{12}>0$; this follows immediately if we
 take expectations to $X_{1:2}X_{2:2}=X_1 X_2$ and $X_{1:2}+X_{2:2}=X_1+X_2$, yielding
 $\Cov(X_{1:2},X_{2:2})=\left(\E X_{2:2}-\E X_1\right)^2>0$, since
 $\E X_1<\E X_{2:2}$ for any non-degenerate $F$.
 At the time of 80's, David Scott and Robert
 Bartoszy\'{n}ski, in connection with a problem in cell division,
 proposed to George R.\ Terrell the conjecture that $\rho_{12}$ is never
 greater then one-half. This was proved true:

 \begin{theo}{\rm (Terrell, 1983).
 \label{theo.Terrell}
 If $F$ is a continuous distribution with
 finite variance then $\rho_{12}\leq 1/2$, with equality
 if and only if $F$ is a rectangular distribution.
 }
 \end{theo}

 \noindent
 Terrell's result is based on the development of a function in
 a series of Legendre polynomials, it is quite complicated,
 and imposes the unnecessary restriction that $F$ is continuous.
 Sz\'{e}kely and M\'{o}ri (1985)
 interpreted Terrell's result as a maximal correlation problem --
 see Gebelein (1941), R\'{e}nyi (1959) -- that is,
 \[
 R(X_{1:2},X_{2:2}):=\sup_{g_1,g_2} \rho(g_1(X_{1:2}),g_2(X_{2:2}))
 \]
 where the supremum is taken over
 non-constant functions $g_1\in L^2(X_{1:2})$, $g_2\in L^2(X_{2:2})$.
 In this way, Sz\'{e}kely and M\'{o}ri improved Terrell's
 result in four directions.
 First, they removed the restriction that $F$ is continuous;
 second, they simplified Terrell's proof;
 third, they showed that the correlation coefficient of
 $g_1(X_{1:2})$ and $g_2(X_{2:2})$ is less than $1/2$,
 for any distribution and any
 functions $g_1(X_{1:2})$ and $g_2(X_{2:2})$;
 fourth, and most important,
 they extended these results to any sample size $n\geq 2$,
 obtaining the inequality (denote by $X_{1:n}\leq \cdots\leq X_{n:n}$
 the order statistics from a random sample of size $n$ from $F$)
 \be
 \label{SzM1985}
 \rho (X_{i:n},X_{j:n}) \leq \sqrt{\frac{i(n+1-j)}{j(n+1-i)}}, \ \ \
 1\leq i<j\leq n, \ \ n\geq 2,
 \ee
 which is valid for all non-degenerate
 distribution functions $F$ with $\Var X_{i:n}+\Var X_{j:n}<\infty$.
 The equality in (\ref{SzM1985}), for a single value of $(i,j,n)$,
 characterizes the rectangular distributions.

 Subsequently, Nevzorov (1992) obtained a similar inequality for upper records,
 in which the equality characterizes the location-scale family
 of the standard exponential random variable.
 Later on, L\'{o}pez-Bl\'{a}zquez and Casta\~{n}o-Mart\'{i}nez (2006)
 proved that the TSM inequality (\ref{SzM1985}) is also valid
 when the order statistics are based on a without-replacement sample of size $n$,
 taken from a finite population ${\cal P}=\{x_1,\ldots,x_N\}$ with $N>n$
 distinct elements.

 It is shown in Papadatos and Xifara (2013) that all
 presenting results are based, essentially, in the following
 polynomial regression property (PRP):
 \be
 \label{Polynomial.reggression}
 \E(X^k|Y)=A_k Y^k +P_{k-1}(Y), \ \ \ \E(Y^k|X)=B_k X^k+Q_{k-1}(X), \ \ k=1,2,\ldots,
 \ee
 where $P_{k-1}$ and $Q_{k-1}$ are polynomials of degree at most $k-1$.
 If a random pair $(X,Y)$ satisfies (\ref{Polynomial.reggression})
 then, under mild conditions, its maximal correlation, $R(X,Y)$,
 equals to $\sup_{k\geq 1}\sqrt{A_k B_k}$, so one simply has to
 calculate the principal coefficients $A_k$ and $B_k$ in
 (\ref{Polynomial.reggression}), and choose the value of $k=k_0$
 that maximizes the products $A_kB_k$. If this value of $k_0$
 is unique, then the equality in the inequality
 \[
 \rho(g_1(X),g_2(Y))\leq R(X,Y)
 \]
 is attained if and only if $\Pr[g_1(X)=c_1+\lambda_1 \phi_{k_0}(X)]=1$,
 $\Pr[g_2(Y)=c_2+\lambda_2 \psi_{k_0}(Y)]=1$, and
 $\lambda_1\lambda_2\sign(A_{k_0})>0$,
 where $\{\phi_k\}_{k=0}^\infty$ is the complete, orthonormal
 polynomial system
 in $L^{2}(X)$ (with the convention that each $\phi_k$ has positive
 principal coefficient)
 and $\{\psi_k\}_{k=0}^\infty$ the corresponding system
 in $L^{2}(Y)$.
 Let us denote by $U_{1:n}<\cdots<U_{n:n}$ the order statistics 
 either
 from uniform in the interval $(0,1)$, or from the discrete
 uniform in $\{1,\ldots,N\}$ (the latter in the without-replacement case).
 It is a simple exercise to verify
 that the pair $(U_{i:n},U_{j:n})$ posses the PRP
 (\ref{Polynomial.reggression}),
 and it is easy to calculate $A_k$ and $B_k$
 (in the without-replacement case this is slightly more complicated);
 also, it is plain to check the PRP for upper records $(W_n,W_m)$
 ($n<m$) from the standard exponential, and to obtain the
 constants $A_k$ and $B_k$ in (\ref{Polynomial.reggression}).
 Then, the PRP method yields all the presenting results
 at once and, in the discrete case, it provides the exact
 characterization of those populations ${\cal P}$
 that attain the upper bound in the TSM inequality (\ref{SzM1985}),
 even if ties are allowed;
 see Section 3 in Papadatos and Xifara (2013).
 In all of the forgoing results, the sequence
 $A_kB_k$ is uniquely maximized by its first term,
 hence, $R=\sqrt{A_1 B_1}$, and the equality is attained by
 linear functions $g_1$, $g_2$; equivalently,
 $R(X,Y)=|\rho(X,Y)|$. While the equality $|\rho|=R$
 often appears to problems regarding maximal correlation
 under PRP, it is not always true; see Papadatos (2014).

 In the present note we shall prove the following discrete 
 analogue of Theorem \ref{theo.Terrell}. First, we provide
 a useful definition, cf.\ Balakrishnan {\it et al.}\ (2003).
  
 \begin{defi}{\rm
 \label{def.ties.finite}
 (Uniform distribution on
 ${\cal P}=\{x_1,\ldots,x_N\}$ in the possible presence of ties). 
 Let $k_1$ of the $x_i$'s be equal to $y_1$, $k_2$ of the $x_i$'s
 be equal to $y_2$, and so on, where $k_1+\cdots+k_m=N$ and,
 with no loss of generality, $y_1<\cdots<y_m$.
 We say that $X$ is uniform on ${\cal P}=\{x_1,\ldots,x_N\}$
 if $\Pr(X=y_j)=k_j/N$, $j=1,\ldots,m$.
 }
 \end{defi}

 The main result is the following.
 
 \begin{theo}{\rm
 \label{theo.rho}
 Let $X_1$, $X_2$ be independent random variables with uniform
 distribution on ${\cal P}=\{x_1,\ldots,x_N\}$, where $N\geq 2$ and
 $x_1\leq \cdots \leq x_N$ with $x_1<x_N$.
 Then, under the notation of Theorem
 \ref{theo.Terrell},
 \be
 \label{rho.new}
 \rho_{12}\leq \frac{1-N^{-2}}{2+N^{-2}},
 \ee
 and the equality is attained if and only if
 ${\cal P}$  is a discrete lattice, that is,
 $x_{i+1}-x_i=\lambda>0$
 (constant), $i=1,\ldots,N-1$.
 }
 \end{theo}

 \noindent
 A straightforward computation shows that the upper bound
 in (\ref{rho.new}) equals to $\rho_{12}$
 when ${\cal P}={\cal P}_0:=\{1,\ldots,N\}$. If $U_{1:2}\leq U_{2:2}$
 are the order statistics from ${\cal P}_0$, it can be checked
 that the random pair $(U_{1:2},U_{2:2})$ does not posses
 the PRP (\ref{Polynomial.reggression}); in particular,
 $U_{1:2}$ does not have linear regression on $U_{2:2}$. Indeed,
 for $x,y\in{\cal P}_0$ with $x\leq y$,
 $\Pr(U_{1:2}=x|U_{2:2}=y)=2/(2y-1)$ if $x<y$,
 and $\Pr(U_{1:2}=y|U_{2:2}=y)=1/(2y-1)$. Hence, one finds
 $\E (U_{1:2}|U_{2:2}=y)=y^2/(2y-1)$, and this regression
 (with $y$ restricted to ${\cal P}_0$) is, clearly, nonlinear, unless $N=2$.
 Therefore, (\ref{Polynomial.reggression}) fails, and this
 verifies the essentially different nature of the present problem, compared
 to the preceding ones. Notice that the case $N=2$ is trivial, since
 $\rho_{12}\equiv 1/3$ for every choice of $x_1$, $x_2$, with $x_1<x_2$
 (the coefficient of correlation is location-scale invariant).

 \section{Proof of Theorem \ref{theo.rho}}
 \label{sec.2}
 We shall apply the ordinary method of Terrell (1983) to the case
 where $U$ follows a discrete uniform in ${\cal P}_0=\{1,\ldots,N\}$.
 In this case, the orthonormal polynomial system is provided
 by Hahn polynomials, $\{\psi_k\}_{k=0}^{N-1}$,
 where
 \[
 \psi_k(x)=N^{1/2}{N+k\choose 2k+1}^{-1/2}{2k\choose k}^{-1/2}
 \sum_{j=0}^k (-1)^{k-j} {k+j \choose j}
 {N-1-j \choose k-j} {x-1\choose j}.
 \]
 The following two properties will be used in 
 the sequel; cf.\  L\'{o}pez-Bl\'{a}zquez and 
 Casta\~{n}o-Mart\'{i}nez (2006).
 First,
 \[
 \psi_k(x)=A_k x^k + B_{k} x^{k-1}+\cdots,
 \]
 where
 \be
 \label{Hahn.coefficients}
  A_k=\frac{\sqrt{N}}{k!}{2k\choose k}^{1/2}
 {N+k\choose 2k+1}^{-1/2},
 \ \ \ \ \
 B_k=-\frac{k(N+1)\sqrt{N}}{k!}
 {2k-1 \choose k} {N+k\choose 2k+1}^{-1/2} {2k\choose k}^{-1/2}.
 \ee
 Second, the polynomials $\psi_k$ are orthonormal in
 $L^2(U)\equiv\{g:{\cal P}_0\to\RR\}$, that is,
 \[
 \E \psi_{k}(U)\psi_m(U)=\frac{1}{N}\sum_{x=1}^N \psi_{k}(x)\psi_m(x)
 =\delta_{k,m}=\left\{\begin{array}{ccc}
 1, &\mbox{if}& k=m, \\
 0, &\mbox{if}& k\neq m,
 \end{array}\right.
 \]
 and this implies that any function $g:{\cal P}_0\to\RR$ has a unique
 development in Hahn polynomial series, namely,
 \be
 \label{Hahn.series}
 g(x)=\sum_{k=0}^{N-1} \delta_k \psi_k(x), \ \ \ \mbox{where}
 \ \ \delta_k=\E \psi_k(U) g(U) \ \
 (k=0,\ldots,N-1).
 \ee
 The constants $\{\delta_k\}_{k=0}^{N-1}$ are the Fourier coefficients of $g$.

 An arbitrary ${\cal P}=\{x_1,\ldots,x_N\}$ (with $x_1\leq \cdots\leq x_N$)
 is the image of some nondecreasing $g:{\cal P}_0\to\RR$,
 namely, $g(i)=x_i$, $i=1,\ldots,N$. In other words, if $X$ is uniformly
 distributed over ${\cal P}$ (in the sense of Definition \ref{def.ties.finite}) 
 and $U$ is uniformly distributed over
 ${\cal P}_0$ then $X$ and $g(U)$ have the same distribution.
 Therefore, $(g(U_1),g(U_2))$ and $(X_1,X_2)$ have the same distribution,
 where $X_1,X_2$ are independent copies of $X$, and $U_1,U_2$ are
 independent copies of $U$.
 Since $g$ is 
 nondecreasing, it holds $(\min\{g(U_1),g(U_2)\},\max\{g(U_1),g(U_2)\})
 =(g(\min\{U_1,U_2\}),g(\max\{U_1,U_2\}))=(g(U_{1:2}),g(U_{2:2}))$.
 This shows that the random pairs $(X_{1:2},X_{2:2})$ and $(g(U_{1:2}),g(U_{2:2}))$
 have the same distribution;
 hence, $\rho_{12}=\rho(g(U_{1:2}),g(U_{2:2}))$.
 Expanding $g$ according to (\ref{Hahn.series}), we will be able to
 express $\rho_{12}$ as a function of $\delta_1,\ldots,\delta_{N-1}$,
 and, next, we shall maximize the resulting function. The calculations,
 below, do not a-priori impose any monotonicity restrictions on $g$;
 the only assumption is that $g$ is non-constant (otherwise, the variance
 of $X$ is zero and $\rho_{12}$ is undefined).

 \begin{lem}{\rm
 (Hahn representation).
 \label{lem.s11.s22.s12}
 Expand an arbitrary $g:{\cal P}_0\to\R$ in a Hahn polynomial series
 as in (\ref{Hahn.series}),
 and let $\sigma_1^2=\Var g(U_{1:2})$,
 $\sigma_2^2=\Var g(U_{2:2})$,
 $\sigma_{12}=\Cov(g(U_{1:2}),g(U_{2:2}))$. Then,
 \begin{eqnarray}
 \label{s11}
 \sigma_1^2
 &\hspace*{-1ex}=\hspace*{-1ex}&
 \sum_{k=1}^{N-1} \delta_k^2 -\frac{\delta_1^2}{3}\Big(1-\frac{1}{N^2}\Big)
 -\frac{2}{N}\sum_{k=1}^{N-2}\lambda_k \delta_k\delta_{k+1}
 \\
 \label{s22}
 \sigma_2^2
 &\hspace*{-1ex}=\hspace*{-1ex}&
 \sum_{k=1}^{N-1} \delta_k^2 -\frac{\delta_1^2}{3}\Big(1-\frac{1}{N^2}\Big)
 +\frac{2}{N}\sum_{k=1}^{N-2}\lambda_k \delta_k\delta_{k+1}
 \\
 \label{s12}
 \sigma_{12}
 &\hspace*{-1ex}=\hspace*{-1ex}&
 \frac{\delta_1^2}{3}\Big(1-\frac{1}{N^2}\Big),
 \end{eqnarray}
 where, in (\ref{s11}) and (\ref{s22}),
 \be
 \label{lambda.k}
 \lambda_{k}=(k+1)\sqrt{\frac{N^2-(k+1)^2}{(2k+1)(2k+3)}}, \ \ \ k=1,\ldots,N-1.
 \ee
 }
 \end{lem}
 \begin{pr}{Proof}
 If $\delta_1=\cdots=\delta_{N-1}=0$ then $g\equiv\delta_0$ (constant)
 and the result is trivial, since $\sigma_1^2=\sigma_2^2=\sigma_{12}=0$.
 Assume that $\beta:=\delta_1^2+\cdots+\delta_{N-1}^2>0$, and set
 $h(x):=(g(x)-\delta_0)/\sqrt{\beta}$, so that
 \[
 h(x)=\sum_{k=1}^{N-1} d_k \psi_k(x), \ \  \mbox{ where }
 \ d_k=\delta_k/\sqrt{\beta} \ \
 (k=1,\ldots, N-1),
 \]
 are the Fourier coefficients of $h$. By construction,
 $\sum_{k=1}^{N-1}d_k^2=1$,
 $\sigma_{1}^2=\beta \Var h(U_{1:2})$,
 $\sigma_{2}^2=\beta \Var h(U_{2:2})$,
 $\sigma_{12}=\beta \Cov(h(U_{1:2}),h(U_{2:2}))$.
 The orthonormality of $\psi_k$ shows that
 the random variable $Y:=h(U)$ is standardized; indeed,
 \[
 \E Y =\sum_{k=1}^{N-1} d_k \E \psi_k(U)=0, \ \ \ \
 \E Y^2=\sum_{k=1}^{N-1} \sum_{m=1}^{N-1}d_k d_m
 \E \psi_k(U)\psi_m(U)=\sum_{k=1}^{N-1}d_k^2\E \psi_k(U)^2=1,
 \]
 where, the first equality is satisfied because $\psi_k$
 ($k\geq 1$) is orthogonal to $\psi_0\equiv 1$, and the second one
 follows  from the orthonormality of $\psi_k$ and the definition
 of $d_k$.
 We now write $Y_1=h(U_1)$, $Y_2=h(U_2)$, $Z_1=h(U_{1:2})$,
 $Z_2=h(U_{2:2})$, with $U_1$, $U_2$ being independent uniform
 from ${\cal P}_0$; notice that the inequality $Z_1\leq Z_2$
 may fail, since $h$ has not been assumed monotonic.
 According to the argument above, $Y_1,Y_2$ are independent
 and standardized. Taking expectations
 to the obvious identities $Z_1+Z_2=Y_1+Y_2$ and
 $Z_1 Z_2=Y_1 Y_2$ we see that $\E Z_1=-\E Z_2$ and $\E Z_1 Z_2=0$.
 Hence, $\Cov(Z_1,Z_2)=(\E Z_2)^2$; equivalently,
 $\sigma_{12}=\beta (\E Z_2)^2$.
 The relation $2=\Var(Y_1+Y_2)=\Var(Z_1+Z_2)=\Var Z_1+\Var Z_2+2 (\E Z_2)^2$
 shows that $\Var Z_1=2-\Var Z_2-2(\E Z_2)^2=2-(\E Z_2)^2-\E Z_2^2$ and
 $\sigma_1^2=\beta\big(2-(\E Z_2)^2-\E Z_2^2\big)$.
 Since $\Var Z_2=\E Z_2^2-(\E Z_2)^2$ implies
 $\sigma_2^2=\beta\big(\E Z_2^2-(\E Z_2)^2\big)$, all
 quantities $\sigma_1^2$, $\sigma_2^2$, $\sigma_{12}$, 
 (and $\rho_{12}$),
 are expressed in terms of $\E Z_2=\E h(U_{2:2})$
 and $\E Z_2^2=\E h(U_{2:2})^2$.

 The probability mass function of $U_{2:2}$ is given by
 $\Pr(U_{2:2}=j)=(2j-1)/N^2$, $j=1,\ldots,N$, and therefore,
 \[
 \E Z_2=\frac{1}{N^2}\sum_{j=1}^{N} (2j-1)h(j)=\frac{1}{N}\E (2U-1)h(U)
 =\frac{2}{N}\E U h(U),
 \]
 where the last equality follows from $\E h(U)=\E Y=0$. Since
 each $\psi_k(U)$ ($k\geq 2$) is orthogonal to any polynomial of degree
 at most one, hence to $U$, we have
 \[
 \E U h(U) =d_1 \E U\psi_1(U)
 +\sum_{k=2}^{N-1}d_k\E U\psi_k(U)= d_1\E U\psi_1(U)
 =\frac{d_1}{A_1}\E \psi_1(U) (A_1 U+B_1)=\frac{d_1}{A_1},
 \]
 because $A_1 U+B_1=\psi_1(U)$ and $\psi_1$ is orthogonal to constants.
 Substituting the value of $A_1$
 from (\ref{Hahn.coefficients}) we obtain
 $\E Z_2=2d_1/(N A_1)=(1-N^{-2})^{1/2}d_1/\sqrt{3}$, and the relations
 $\sigma_{12}=\beta (\E Z_2)^2$ and $d_1=\delta_1/\sqrt{\beta}$
 imply
 (\ref{s12}).

 The computation of the second moment of $Z_2$ is more involved. Write
 \be
 \label{aux.1}
 \E Z_2^2=\frac{1}{N^2}\sum_{j=1}^N (2j-1) h(j)^2=
 \frac{1}{N}\E (2 U-1) h(U)^2=\frac{2}{N}\E U h(U)^2-\frac{1}{N},
 \ee
 because $\E Y^2=\E h(U)^2=1$. It remains to compute
 \[
 \E U h(U)^2=\sum_{k,m=1}^{N-1}d_k d_m
 \E U \psi_k(U) \psi_m(U)
 =
 \sum_{k=1}^{N-1} d_k^2 \E U \psi_k(U)^2 +
 2\sum_{k=1}^{N-2} d_k d_{k+1} \E U \psi_k(U)\psi_{k+1}(U);
 \]
 all other terms (with $m>k$) vanish, due to the orthogonality of $\psi_m(U)$
 to the polynomial $U\psi_k(U)$ (of degree $k+1$) when $m\geq k+2$.
 We proceed to compute $\E U \psi_k(U)^2$ and $\E U \psi_k(U) \psi_{k+1}(U)$.
 Write $\psi_k(x)=A_k x^k +B_k x^{k-1}+P_{k-2}(x)$, where
 $P_{k-2}$ is a polynomial of degree at most $k-2$  ($P_{-1}\equiv 0$).
 Then,
 \[
 \E U\psi_k(U)^2=\E U \psi_k(U)(A_k U^k+B_k U^{k-1}+P_{k-2}(U))
 =A_k \E  U^{k+1} \psi_k(U)+B_k \E U^k\psi_k(U),
 \]
 because $\psi_k(U)$ is orthogonal to the polynomial $U P_{k-2}(U)$.
 Next, we compute
 \begin{eqnarray*}
 \E  U^{k+1} \psi_k(U)
 &\hspace*{-1ex}=\hspace*{-1ex}&
 \frac{1}{A_{k+1}}
 \E \psi_k(U) \Big(A_{k+1}U^{k+1}+B_{k+1} U^{k}+P_{k-1}(U)-B_{k+1}U^k-P_{k-1}(U)\Big)
 \\
 &\hspace*{-1ex}=\hspace*{-1ex}&
 \frac{1}{A_{k+1}}
 \E \psi_k(U) \Big(\psi_{k+1}(U)-B_{k+1}U^k-P_{k-1}(U)\Big)
 \\
 &\hspace*{-1ex}=\hspace*{-1ex}&
 -\frac{B_{k+1}}{A_{k+1}}
 \E U^k \psi_k(U),
 \end{eqnarray*}
 where the last equality follows from the orthogonality of $\psi_k$
 to both $\psi_{k+1}$ and $P_{k-1}$. Therefore,
 \[
 \E U\psi_k(U)^2=A_k \E  U^{k+1} \psi_k(U)+ B_k \E U^k\psi_k(U)=
  \left(B_k-\frac{A_k B_{k+1}}{A_{k+1}}\right) \E U^k \psi_k(U).
 \]
 The calculation of $\E U^k \psi_k(U)$ is easy:
 \[
 \E U^k \psi_k(U)=\frac{1}{A_k}
 \E\psi_k(U) (A_k U^k+B_{k-1} U^{k-1}+\cdots)
 =\frac{1}{A_k}\E \psi_k(U)^2=\frac{1}{A_k}.
 \]
 Hence, using the fact that $B_k/A_k=-k(N+1)/2$,
 see (\ref{Hahn.coefficients}), we obtain
 \[
  \E U\psi_k(U)^2=\frac{B_k}{A_k}-\frac{B_{k+1}}{A_{k+1}}=\frac{N+1}{2};
 \]
 surprisingly, this expectation is independent of $k$.
 Proceeding similarly, and in view of (\ref{Hahn.coefficients}), we  compute
 \[
 \E U\psi_k(U)\psi_{k+1}(U)=\E U (A_k U^k+B_k U^{k-1}+\cdots) \psi_{k+1}(U)
 =A_k \E U^{k+1}\psi_{k+1}(U)=\frac{A_k}{A_{k+1}}=\frac{\lambda_k}{2},
 \]
 with $\lambda_k$ as in (\ref{lambda.k}).
 Combining the preceding formulae (recall that $\sum_{k=1}^{N-1}d_k^2=1$)
 we obtain
 \[
 \E U h(U)^2
 =
 \frac{N+1}{2}+
  \sum_{k=1}^{N-2}\lambda_k
  d_k d_{k+1}.
 \]
 Then, from (\ref{aux.1}),
 $\E Z_2^2 = 1  + (2/N) \sum_{k=1}^{N-2}  \lambda_k  d_k d_{k+1}$.
 Moreover, since  $(\E Z_2)^2=(1-N^{-2})d_1^2/3$ we have
 $\Var Z_2=1-(1-N^{-2})d_1^2/3  + (2/N) \sum_{k=1}^{N-2}
 \lambda_k  d_k d_{k+1}$, and
 the relation $\sigma_2^2=\beta\Var Z_2$ yields (\ref{s22}), because
 $\beta d_1^2=\delta_1^2$ and  $\beta d_{k} d_{k+1}=\delta_k \delta_{k+1}$
 (recall  that $\beta=\delta_1^2+\cdots+\delta_{N-1}^2$).
 Similarly, substituting the preceding expressions for $(\E Z_2)^2$ and $\E Z_2^2$
 to the relation $\Var Z_1=2-(\E Z_2)^2-\E Z_2^2$ we obtain
 $\Var Z_1=1-(1-N^{-2})d_1^2/3-(2/N) \sum_{k=1}^{N-2}  \lambda_k  d_k d_{k+1}$,
 and since $\sigma_1^2=\beta\Var Z_1$, (\ref{s11}) follows.
 \end{pr}

 \begin{lem}{\rm
 \label{lem.rational.functions}
 We define the rational functions $R_0(x)\equiv 1$
 and, recurrently,
 \[
 R_k(x):=1-\frac{k^2 (1-k^2/x)}{(4 k^2-1)R_{k-1}(x)},
 \ \ \ \ k=1,2,\ldots.
 \]
 Then, the following inequalities hold true:
 \be
 \label{rational.inequalities}
 \frac{k+1}{2 k+1} < R_k(x)<1, \ \  x>k^2, \ \ k=1,2,\ldots.
 \ee
 }
 \end{lem}
 \begin{pr}{Proof}
 It is easily seen that $R_1(x)=(2+x^{-1})/3$, which is strictly decreasing
 on
 $[1,\infty)$, with $R_1(1)=1$ and $R_1(\infty)=2/3$; hence,
 (\ref{rational.inequalities}) is fulfilled for $k=1$.
 Suppose that (\ref{rational.inequalities}) is true for some $k-1\geq 1$, that
 is, $k/(2k-1)<R_{k-1}(x)<1$, $x>(k-1)^2$;
 equivalently,
 \[
 1<\frac{1}{R_{k-1}(x)}<\frac{2k-1}{k}, \ \ \ x>(k-1)^2.
 \]
 If $x>k^2$, multiplying the preceding inequality
 by the positive quantity $k^2 (1-k^2/x)/(4 k^2-1)$
 we get
 \[
 L(x):=\frac{k^2 (1-k^2/x)}{4 k^2-1}
 <\frac{k^2 (1-k^2/x)}{(4 k^2-1)R_{k-1}(x)}
 =1-R_k(x)<U(x):=\frac{k(1-k^2/x)}{2k+1}, \ \ \ x>k^2.
 \]
 It is easily seen that $L$ is strictly increasing in $x$, $x>k^2$,
 so $L(x)>L(k^2)=0$, and $U$ is also strictly increasing with
 $U(x)<U(\infty)=k/(2k+1)$. This shows that $0<1-R_k(x)<k/(2k+1)$
 for $x>k^2$ and verifies the inductional step,
 completing the proof.
 \end{pr}

 \begin{lem}{\rm
 \label{lem.numbers.ak.bk}
 For $N\geq 2$, there exist nonnegative numbers
 $\{\alpha_k\}_{k=1}^{N-1}$ and $\{\beta_k\}_{k=1}^{N-1}$
 satisfying the following four properties:
 \begin{itemize}

 \item[(i)]
 $\alpha_1=\frac{1}{3}\Big(2+\frac{1}{N^2}\Big)$, \ \ $\beta_{N-1}=0$;

 \item[(ii)]
 $\alpha_k>\beta_{k-1}> 0$, \ \ $k=2,\ldots,N-1$;

 \item[(iii)]
 $\alpha_k+\beta_{k-1}=1$, \ \ $k=2,\ldots,N-1$;

 \item[(iv)]
 $\alpha_k\beta_k=\frac{(k+1)^2}{(2k+1)(2k+3)}\Big(1-\frac{(k+1)^2}{N^2}\Big)$,
 \ \ \ $k=1,\ldots,N-1$.
 \end{itemize}
 }
 \end{lem}
 \begin{pr}{Proof} Define $\alpha_k:=R_k(N^2)$ and $\beta_k:=1-R_{k+1}(N^2)$,
 $k=1,\ldots,N-1$, where the functions $R_k$ are
 defined in Lemma \ref{lem.rational.functions}. Then, (i)
 and (iii) are obvious, and (ii) follows from (\ref{rational.inequalities}),
 since $\alpha_k=R_k(N^2)>(k+1)/(2k+1)>1/2$, because $N^2>k^2$,
 and, by the same reasoning, $\beta_k\in[0,1/2)$ (more precisely,
 $\beta_k>0$ unless $k=N-1$). Finally, (iv) is trivial
 for $k=N-1$ (since $\beta_{N-1}=0$), and for $k=1,\ldots,N-2$,
 \[
 \alpha_k\beta_k=\alpha_k (1-\alpha_{k+1})=
 \alpha_k\frac{(k+1)^2(1-(k+1)^2/N^2)}{[4(k+1)^2-1] \alpha_k}
 =
 \frac{(k+1)^2}{(2k+1)(2k+3)}\Big(1-\frac{(k+1)^2}{N^2}\Big),
 \]
 due to the recurrent relation of Lemma \ref{lem.rational.functions},
 namely,
 \[
 1-R_{k+1}(N^2)
 =\frac{(k+1)^2(1-(k+1)^2/N^2)}{[4(k+1)^2-1] R_k(N^2)}.
 \]
 Hence, the lemma is proved.
 \end{pr}

 \begin{lem}{\rm (Terrell-Hahn representation).
 \label{lem.Terrell.Hann}
 Consider the numbers $\{\alpha_k\}_{k=1}^{N-1}$ and
 $\{\beta_k\}_{k=1}^{N-1}$ as in Lemma \ref{lem.numbers.ak.bk}.
 Then, with the convention $\delta_N:=0$,
 \[
 \sigma_1^2=\sum_{k=1}^{N-1} \left(\sqrt{\alpha_k}\delta_k
 -\sqrt{\beta_k}\delta_{k+1}\right)^2, \ \ \ \
 \sigma_2^2=\sum_{k=1}^{N-1} \left(\sqrt{\alpha_k}\delta_k
 +\sqrt{\beta_k}\delta_{k+1}\right)^2,
 \]
 where $\sigma_1^2=\Var g(U_{1:2})$  and $\sigma_2^2=\Var g(U_{2:2})$
 are as in Lemma \ref{lem.s11.s22.s12}.
 }
 \end{lem}
 \begin{pr}{Proof}
 Let $s_1^2$ be the first sum. Expanding the squares we have
 \[
 s_1^2=\sum_{k=1}^{N-1} \alpha_k \delta_k^2
  +\sum_{k=1}^{N-1} \beta_k \delta_{k+1}^2
  -2\sum_{k=1}^{N-1}\sqrt{\alpha_k\beta_k}\delta_k\delta_{k+1}
  =\alpha_1\delta_1^2
  +\sum_{k=2}^{N-1} (\alpha_k+\beta_{k-1}) \delta_k^2
  -2\sum_{k=1}^{N-2}\sqrt{\alpha_k\beta_k}\delta_k\delta_{k+1},
 \]
 and Lemma \ref{lem.numbers.ak.bk} implies that
 \[
 s_1^2=\frac{1}{3}\Big(2+\frac{1}{N^2}\Big)\delta_1^2
  +\sum_{k=2}^{N-1} \delta_k^2
  -\frac{2}{N}\sum_{k=1}^{N-2}
  \lambda_k \delta_k\delta_{k+1},
 \]
 with $\lambda_k$ as in (\ref{lambda.k}).
 Observing that $(2+N^{-2})/3=1-(1-N^{-2})/3$, $s_1^2$ reduces to
 the expression of $\sigma_1^2$, given in (\ref{s11}).
 The derivation of the expression for $\sigma_2^2$
 is completely similar.
 \end{pr}
 \begin{theo}{\rm
 \label{theo.rho.general}
 For $N\geq 2$ and any nonconstant $g:{\cal P}_0=\{1,\ldots,N\}\to\R$,
 \[
 \rho(g(U_{1:2}),g(U_{2:2}))\leq \frac{1-N^{-2}}{2+N^{-2}},
 \]
 with equality if and only if $g$ is linear.
 \label{theo.genaral}
 }
 \end{theo}
 \begin{pr}{Proof}
 Let $g(U)=\sum_{k=0}^{N-1}\delta_k \psi_k(U)$ be the Hahn representation
 of $g$, with $\delta_1^2+\cdots+\delta_{N-1}^2>0$ (otherwise, $g$ is constant).
 If $\sigma_i^2=\Var g(U_{i:2})$ ($i=1,2$) and
 $\sigma_{12}=\Cov(g(U_{1:2}),g(U_{2:2}))$
 then, from Lemma \ref{lem.Terrell.Hann} (and with the convention $\delta_N=0$),
 \begin{eqnarray*}
 \sigma_1^2 \sigma_2^2
 &\hspace*{-1ex}=\hspace*{-1ex}&
 \Bigg\{\sum_{k=1}^{N-1} \left(\sqrt{\alpha_k}\delta_k
 -\sqrt{\beta_k}\delta_{k+1}\right)^2\Bigg\}
 \Bigg\{\sum_{k=1}^{N-1} \left(\sqrt{\alpha_k}\delta_k
 +\sqrt{\beta_k}\delta_{k+1}\right)^2\Bigg\}
 \\
 &\hspace*{-1ex}\geq\hspace*{-1ex}&
 \Bigg\{\sum_{k=1}^{N-1}
 \left(\sqrt{\alpha_k}\delta_k
 -\sqrt{\beta_k}\delta_{k+1}\right)
 \left(\sqrt{\alpha_k}\delta_k
 +\sqrt{\beta_k}\delta_{k+1}\right)\Bigg\}^2
 \\
 &\hspace*{-1ex}=\hspace*{-1ex}&
 \Bigg\{\sum_{k=1}^{N-1}
 \left(\alpha_k\delta_k^2
 -\beta_k\delta_{k+1}^2\right)\Bigg\}^2
 \\
 &\hspace*{-1ex}=\hspace*{-1ex}&
 \Bigg\{\alpha_1\delta_1^2+\sum_{k=2}^{N-1}
 \left(\alpha_k-\beta_{k-1}\right)\delta_k^2
 \Bigg\}^2
 \\
 &\hspace*{-1ex}\geq\hspace*{-1ex}&
 \alpha_1^2\delta_1^4, \Bigg.
 \end{eqnarray*}
 where the second inequality holds true because $\alpha_k-\beta_{k-1}>0$
 and $\alpha_1>0$,
 see Lemma \ref{lem.numbers.ak.bk}(ii), while the first one is a simple
 application of the well-known Cauchy inequality. On substituting the value
 of $\alpha_1$ from Lemma \ref{lem.numbers.ak.bk}(i), the previous
 argument verifies the inequality
 \[
 \left(\sigma_1^2\sigma_2^2\right)^{1/2}\geq
 \frac{1}{3}\Big(2+\frac{1}{N^2}\Big)\delta_1^2,
 \]
 in which the equality holds if and only if $\delta_k=0$ for all $k\geq 2$, i.e.,
 when $g$ is linear. From (\ref{s12}) of Lemma \ref{lem.s11.s22.s12},
 $\delta_1^2=3\sigma_{12}/(1-N^{-2})$; thus,
 $(\sigma_1^2\sigma_2^2)^{1/2}\geq
 \sigma_{12} (2+N^{-1})/(1-N^{-2})$, and this
 is equivalent to the desired inequality.
 \medskip
 \end{pr}

 Applying Theorem \ref{theo.rho.general} to the particular 
 case where $g$ is
 given by $g(i)=x_i$, $i=1,\ldots,N$ (which is nondecreasing and non-constant),
 we conclude the result of Theorem \ref{theo.rho}.

 \section{Concluding Remarks}
 \label{sec.3}
 \begin{rem}{\rm
 \label{rem.HDG}
 (An improvement to the covariance/expectation bound for $n=2$).
 Since the variance of $X=g(U)$ equals to
 $\sigma^2=\delta_1^2+\cdots+\delta_{N-1}^2\geq \delta_1^2$,
 an immediate consequence of (\ref{s12}) is the Bessel-type
 inequality $\Cov(X_{1:2},X_{2:2})\leq \sigma^2(1-N^{-2})/3$;
 this is similar to the bound
 $\Cov(X_{1:2},X_{2:2})\leq \sigma^2(1-2 N^{-1})/(3(1-N^{-1}))$,
 obtained by Balakrishnan {\it et al.}\ (2003) in the
 without-replacement case. Notice that the equality in any of these
 two bounds is attained only for lattice populations.
 Both bounds converge to $\sigma^2/3$ as $N\to\infty$, and the
 maximizing population (provided $\mu=\E X$ and $\sigma^2$ are held
 fixed), converges weakly to the uniform
 ${\cal U}(\mu-\sigma\sqrt{3},\mu+\sigma\sqrt{3})$  distribution.
 Papathanasiou (1990) proved that
 the inequality $\Cov(X_{1:2},X_{2:2})\leq \sigma^2/3$ holds
 true
 for random samples $X_1,X_2$ from any distribution with
 finite variance $\sigma^2$,
 and the equality characterizes the location-family of
 rectangular distributions with fixed variance.
 Later on, Balakrishnan and Balasubramanian (1993) observed
 that, in view of the obvious identity
 $\Cov(X_{1:2},X_{2:2})=(\E X_{2:2}-\mu)^2$,
 the supremum of this covariance equals to the
 smallest upper bound
 for $\E X_{2:2}-\mu$ in terms of $\sigma$.
 The form of the latter bound, called HDG bound, is long known 
 for $\E X_{n:n}$ (any $n\geq 2$), and it is proved 
 in two classical papers (appeared in consecutive pages
 of the same journal!) by Hartley and David (1954)
 and Gumbel (1954). The particular case ($n=2$) of the HDG bound
 yields the inequality $\E X_{2:2}\leq \mu+\sigma/\sqrt{3}$,
 in which the equality characterizes
 ${\cal U}(\mu-\sigma\sqrt{3},\mu+\sigma\sqrt{3})$.
 Clearly, this is equivalent to Papathanasiou's covariance bound. 
 The improvement
 on the HDG bound for uniform populations (with $N$ elements, in the
 sense of Definition \ref{def.ties.finite}) is an immediate consequence
 of the discrete covariance bound. More precisely,
 \[
 \E X_{2:2}=\mu+(\E X_{2:2}-\mu) =
 \mu+\Big(\Cov(X_{1:2},X_{2:2})\Big)^{1/2}\leq \mu+
 (1-N^{-2})^{1/2} \
 \frac{\sigma}{\sqrt{3}},
 \]
 and the equality characterizes the lattice population with
 mean $\mu$ and variance $\sigma^2$, that is,
 $x_k=\mu+\sigma\sqrt{3}(2k-N-1)/\sqrt{N^2-1}$, $k=1,\ldots,N$.
 }
 \end{rem}

 \begin{rem}{\rm
 \label{rem.rho.functional}
 (Connection to continuous case).
 One may define the functional
 $\rho_{12}(\cdot)$ with domain ${\cal F}:=\{F: F$ is a distribution function
 with finite non-zero variance$\}$ as
 \[
 \rho_{12}(F):=\rho(X_{1:2},X_{2:2}), \ \ \
 \mbox{$X_1$, $X_2$ independent with distribution $F$},
 \]
 where
 $X_{1:2}\leq X_{2:2}$ are the order statistics.
 For $N\geq 2$, we denote by ${\cal F}_N:=\{{\cal P}: {\cal P}=\{x_1,\ldots,x_N\}$ is
 any uniform population with at least two
 distinct elements$\}$ (uniform in the sense of Definition \ref{def.ties.finite}).
 Then it can be verified that
 ${\cal F}_0:=\cup_{N\geq 2}{\cal F}_N$ is dense in ${\cal F}$ with
 respect
 to the weak topology of probability measures, and the functional
 $\rho_{12}$ is continuous.
 Theorem \ref{theo.rho.general} shows that $\sup_{F\in{\cal F}_0}\rho_{12}(F)=1/2$,
 hence, the same is true for $\sup_{F\in{\cal F}}$, and this implies the bound in
 Terrell's result, Theorem  \ref{theo.Terrell}.
 }
 \end{rem}

 \begin{rem}{\rm
 \label{rem.minimization}
 (An equivalent minimization problem).
 To simplify notation, set
 $\gamma_k=2\lambda_k/N$ ($k=1,\ldots,N-1$) with $\lambda_k$
 as in (\ref{lambda.k}),
 and $\lambda=\alpha_1$ as in Lemma \ref{lem.numbers.ak.bk}(i);
 that is,
 \[
 \lambda=\frac{1}{3}\Big(2+\frac{1}{N^2}\Big), \ \ \
 \gamma_k=\frac{2(k+1)}{N}\sqrt{\frac{N^2-(k+1)^2}{(2k+1)(2k+3)}}, \ \
 k=1,\ldots,N-1.
 \]
 According to Lemma \ref{lem.s11.s22.s12}, we have
 \[
 \Big(\frac{1}{\rho_{12}}\Big)^2
 =\frac{9N^4}{(N^2-1)^2} \
 \frac{\Big(\lambda \delta_1^2+\delta_2^2+\cdots+\delta_{N-1}^2\Big)^2-
 \Big(\gamma_1\delta_1\delta_2+\gamma_2\delta_2\delta_3+\cdots+
 \gamma_{N-2}\delta_{N-2}\delta_{N-1}
 \Big)^2}{\delta_1^4},
 \]
 where this expression should be treated as $+\infty$ for $\delta_1=0$.
 The main difficulty in the proof of Theorem \ref{theo.rho.general} 
 was to minimize the above ratio w.r.t.\
 $(\delta_1,\ldots,\delta_{N-1})\in\RR^{N-1}\setminus\{\bbb{0}\}$,
 and the main result 
 showed that the minimizing points are coincide
 with the axis $\{(\delta_1,0,\ldots,0), \delta_1\neq 0\}$.
 The ratio is scale invariant, so, dividing each $\delta_k$
 by $\delta_1\neq 0$ (since $\delta_1=0$ makes
 the ratio infinite, and, certainly, cannot be a minimizer)
 we can assume that $\delta_1=1$. Then, with a slight change
 of notation (remove the constant scalar and set $n=N-1$, $x_k=\delta_k$),
 we define 
 \[
 f(x_2,\ldots,x_n):=
 \Big(\lambda +x_2^2+\cdots+x_n^2\Big)^2-
 \Big(\gamma_1 x_2+\gamma_2 x_2 x_3 +\cdots+\gamma_{n-1}x_{n-1}x_n\Big)^2.
 \]
 Consequently, it remains to minimize the (four degree polynomial)
 function $f$ and to prove that its global minimum is uniquely
 attained at $\bbb{0}$,
 so that $f\geq f(\bbb{0})=\lambda^2$. Despite the fact that
 this (unrestricted) minimization problem looks like a simple exercise,
 this is not so; the minimization points (if exist) depend on the parameters
 $\lambda$, $\gamma_1,\ldots,\gamma_{n-1}$.
 To see this, set $n=3$, $\lambda_0=11/6$ 
 (since $\lambda=\lambda_0$ when $n=3$), 
 and write $f(x_2,x_3;\gamma_1,\gamma_2)
 =(\lambda_0+x_2^2+x_3^2)^2-\big(\gamma_1 x_2+\gamma_2 x_2 x_3\big)^2$.
 It can be checked
 that the function $f(x_2,x_3;4/3,2/3)$,
 though positive, is globally minimized at the points $(x_2,x_3)$ with
 \[
  x_2=\pm \Big\{66 \sqrt{6}-81 \Big\}^{1/2}\Big/16,
  \ \ \
  x_3=  (3 \sqrt{6}-5)\Big/16,
 \]
 not at $\bbb{0}$. Also, $f(\cdot\ ;1,2)$
 is not bounded below; the minimum value of
 $f(\cdot \ ;1,8/5)$ is negative;
 $f(\cdot \ ;1,1)$ attains its global
 minimum at $\bbb{0}$. These examples (in the simplest case $n=3$)
 indicate that
 the multi-parameter case is too complicated, to enable an
 explicit solution. In the contrary, the proposed
 Terrell-type method of proof succeed in obtaining
 Theorems \ref{theo.rho}, \ref{theo.rho.general}.
 Motivated from these examples, a natural question of mathematical nature 
 would ask for (necessary and) sufficient conditions
 on the arbitrary positive parameters $\lambda$, $\gamma_1,\ldots,\gamma_{n-1}$,
 guaranteing  that the minimum of $f(x_2,\ldots,x_n)$, above, is uniquely attained 
 at ${\bbb 0}$.
 }
 \end{rem}

 \begin{rem}{\rm
 \label{rem.max.cor}
 (Connection to maximal correlation).
 It is important to make clear that the assertion of Theorem
 \ref{theo.rho.general}
 cannot be
 reduced to a maximal correlation problem (unless $N=2$),
 in contrast to all of the existing results
 mentioned in Section \ref{sec.1}. To see this,
 it is useful to recall
 R\'{e}nyi's (1959) characterization
 of the maximal correlation of $(X,Y)$, namely,
 \[
 R(X,Y)^2=\mbox{$\sup_f$} \E \Big[\big(\E f(X)|Y\big)^2\Big],
 \]
 where the supremum is taken over $f$ with $\E f(X)=0$,
 $\E f(X)^2=1$. Applying R\'{e}nyi's result to
 $(X,Y)=(U_{1:2},U_{2:2})$, and writing $f(k)=x_k$, $k=1,\ldots,N$,
 it is seen that
 \[
 R^2:=R(U_{1:2},U_{2:2})^2=\frac{1}{N^2}\max \sum_{j=1}^N
 \frac{1}{2j-1} \left\{2(x_1+\cdots+x_{j-1})+x_j\right\}^2,
 \]
 where the maximum is taken
 over $x_1,\ldots,x_N$ satisfying
 \[
 \sum_{k=1}^N (2N-2k+1) x_k=0,  \ \ \ \ \ \sum_{k=1}^N (2N-2k+1) x_k^2=N^2.
 \]
 For $N=3$ one finds a maximization point explicitly, namely,
 \[
 x_1^*=-\sqrt{\frac{2}{5}+\frac{13}{10\sqrt{19}}},
 \ \ \ \ \
 x_2^*=\sqrt{1-\frac{2}{\sqrt{19}}},
 \ \ \ \ \
 x_3^*=\sqrt{4-\frac{1}{2\sqrt{19}}},
 \]
 yielding $R=(2+\sqrt{19})/15>(1-N^{-2})/(2+N^{-2})=8/19$.
 Setting $f(k)=x_k^*$ and $g(Y)=(1/R)\E(f(X)|Y)$,
 that is, $g(k)=y_k^{*}$ with
 \[
 y_1^*=-\sqrt{4-\frac{1}{2\sqrt{19}}},
 \ \ \ \ \
 y_2^*=-\sqrt{1-\frac{2}{\sqrt{19}}},
 \ \ \ \ \
 y_3^*=\sqrt{\frac{2}{5}+\frac{13}{10\sqrt{19}}},
 \]
 if is easy to check that the pair $(f,g)$ is maximally correlated:
 $\E f(U_{1:2})=\E g(U_{2:2})=0$, $\E f(U_{1:2})^2=\E g(U_{2:2})^2=1$, and $\rho(f(U_{1:2}),g(U_{2:2})) =\E f(U_{1:2})g(U_{2:2})=(2+\sqrt{19})/15$.
 
 The situation is similar for all $N$, namely,
 it can be shown that $\rho_N:=(1-N^{-2})/(2+N^{-2})<R_N<1/2$, $N\geq 3$, 
 where $R_N$ is the maximal correlation of an ordered pair from the uniform
 population, ${\cal P}_0$, of size $N$. 
 Clearly, for $N\geq 4$ it would be quite surprising if $R_N$
 could be expressed in a closed form.
 Since the maximal correlation is attainable in our case (the corresponding
 operator is compact; see R\'{e}nyi's, 1959),
 the TSF inequality (\ref{SzM1985}) implies that $R_N< 1/2$.
 Certainly, 
 $R_N\geq \rho_N$. In order to prove that
 $R_N>\rho_N$, it suffices to find a pair of functions with
 coefficient of correlation greater than $\rho_N$.
 Testing small perturbations of linear
 functions, we were led to define $f_0(x):=x-3x^2/N^3$, $g_0(y):=y+3y^2/N^3$
 (both strictly increasing).
 Set $\rho_0=\rho(f_0(U_{1:2}),g_0(U_{2:2}))$,
 $\sigma_1^2=\Var f_0(U_{1:2})$,
 $\sigma_2^2=\Var g_0(U_{2:2})$.
 For $N\geq 3$ we computed (with the help
 of Mathematica)
 \begin{eqnarray*}
 \frac{60^2 \sigma_1^2\sigma_2^2 N^{14} (2 N^2+1)^2}{(N^2 - 1)^2(N^2 - 4)}
 \left( \rho_0^2-\rho_N^2\right)
 &\hspace*{-1ex}=\hspace*{-1ex}&
 -666 - 1332 N + 846 N^2 + 2952 N^3 + 1077 N^4
 \\
 &\hspace*{-1ex} \hspace*{-1ex}&
 - 828 N^5 - 934 N^6 - 792 N^7 - 485 N^8 + 34 N^{10} + 20 N^{12}.
 \end{eqnarray*}
 The substitution  $N=n+3$  ($n\geq 0$)
 shows that the right-hand side,
 $w(n)$, is strictly positive:
 \begin{eqnarray*}
 w(n)
 &\hspace*{-1ex}=\hspace*{-1ex}&
 7010100 + 35183016 n + 72768816 n^2 + 86119956 n^3 + 66523137 n^4 +
 35823456 n^5  \\
 &\hspace*{-1ex} \hspace*{-1ex}&
 + 13910474 n^6 + 3946848 n^7 + 815185 n^8 +
 119820 n^9 + 11914 n^{10} + 720 n^{11} + 20 n^{12}.
 \end{eqnarray*}
 Hence, $R_N\geq \rho_0>\rho_N$ for every $N\geq 3$.
 }
 \bigskip
 \end{rem}

 The main result of the present note, although quite specialized,
 gives rise to a number of questions that, at least to
 author's view, are of some interest.
 Let $U_{1:n}\leq\cdots\leq U_{n:n}$ be the ordered sample
 from the uniform population ${\cal P}_0=\{1,\ldots,N\}$.
 For $1\leq i<j\leq n$ we define
 \[
 R_N=R_N(\mbox{$i$,$j$:$n$})
 :=\mbox{$\sup_{f,g}$} \ \rho(f(U_{i:n}),g(U_{j:n})),  \ \ \ \ \ \
 R'_N=R'_N(\mbox{$i$,$j$:$n$})
 :=\mbox{$\sup_{g}$} \ \rho(g(U_{i:n}),g(U_{j:n})),
 \]
 where the suprema are taken over non-constant $f$, $g$
 (with domain ${\cal P}_0$). The first one is the well-known maximal
 correlation; the second appeared
 in Theorem \ref{theo.rho.general}
 for $i=1$, $j=2$, $n=2$. Moreover, define
 \[
 \rho_N=\rho_N(\mbox{$i$,$j$:$n$}):=\rho(U_{i:n},U_{j:n}),  \ \ \ \ \ \
 R''_N=R''_N(\mbox{$i$,$j$:$n$}):=\mbox{$\sup_{g\nearrow}$} 
 \ \rho(g(U_{i:n}),g(U_{j:n})),
 \]
 where the supremum is taken over nondecreasing non-constant $g$.
 It is obvious that for all $i$, $j$, $n$, $N$,
 $0<\rho_N\leq R_N''\leq R_N'\leq
 R_N<[i(n+1-j)]^{1/2}[j(n+1-i)]^{-1/2}$; see (\ref{SzM1985}).
 Our preliminary calculations indicate the relation
 $\rho_N=R_N''=R_N'<R_N$ for
 all $i$, $j$, $n$, $N$, e.g., 
 $R_N'(\mbox{$1$,$3$:$3$})=(3-7 N^{-2})/(9-N^{-2})$
 (checked for small $N$).
 If the foregoing relation holds true for all choices
 of $i$, $j$, $n$, $N$,
 it will verify that the population ${\cal P}_0$ dominates the
 coefficient of correlation of order statistics from any
 uniform population ${\cal P}=\{x_1,\ldots,x_N\}$, as in
 the continuous case (as well as the without-replacement one).
 An extra difficulty here is due to the lack of a simple formula for 
 $\rho_N(\mbox{$i$,$j$:$n$})$.

 Another natural (but more complicated) extension would be to relax 
 the uniformity assumption. 
 One possible way to do this is to endow the population
 with a fixed probability vector 
 ${\bbb p}=(p_1,\ldots, p_N)$
 with $p_k>0$, $k=1,\ldots,N$.
 Then ${\cal P}_0({\bbb p})=\{1,\ldots,N\}$  is as in the uniform case,
 but, now, the point $k$ has probability $p_k$.
 The problem is to obtain the corresponding results 
 (in an obvious notation) for some of the quantities
 $\rho_N(\mbox{$i$,$j$:$n$};{\bbb p})$, 
 $R''_N(\mbox{$i$,$j$:$n$};{\bbb p})$,
 $R'_N(\mbox{$i$,$j$:$n$};{\bbb p}_N)$, 
 $R_N(\mbox{$i$,$j$:$n$};{\bbb p})$.
 The value of $R''_N$ (and $R_N'$) is of special interest,
 since it provides the least upper bound of the
 correlation of order statistics among all populations
 ${\cal P}({\bbb p})$ endowed with $\bbb p$, i.e., 
 containing $N$ points $x_1\leq \cdots\leq x_N$, $x_1<x_N$, and 
 the probability of $x_k$ is $p_k$.
 Here, in analogy to Definition
 \ref{def.ties.finite}, if ties appear in the population,
 e.g., if $x_{k-1}<x_k=x_{k+1}=\cdots=x_{s}<x_{s+1}$, then
 the probability of $y=x_{k}$
 is defined to be equal to $p_k+\cdots+p_s$.
 In this setup,
 it would be of interest to characterize the population that
 maximizes the correlation of, e.g., the extreme order statistics,
 $\rho\big(X_{1:n},X_{n:n};{\bbb p}\big)$. 
 Clearly, the actual value
 of $R''_N$ depends on ${\bbb p}$, and even the form
 of $R'_3(\mbox{$1$,$2$:$2$};{\bbb p})$ is quite complicated.
 To see this, consider the (binomial-type) vectors
 ${\bbb p}=(1/4,1/2,1/4)$ and ${\bbb p}=(1/16,3/8,9/16)$.
 For the first case we calculate $\rho_3=R'_3=R''_3=9/23$, while
 for the second case one finds
 $\rho_3 =169 \sqrt{3}/\sqrt{655027}=
 0.361674...\hspace*{-.5ex}<0.362354...\hspace*{-.5ex}=R'_3
 =R''_3$,
 where $R_3'$ is the unique positive 
 root of an even, 6th-degree, 
 polynomial, with (large) integer coefficients.
 If all $p_k$ are rational numbers, $p_k=\alpha_k/\beta_k$,
 we can write $p_k=c_k/M$ (in lowest terms). Then,
 Theorem \ref{theo.rho.general} (applied with ties)
 yields the inequality
 $R_N'(\mbox{$1$,$2$:$2$};{\bbb p})\leq 
 R'_M(\mbox{$1$,$2$:$2$})=(1-M^{-2})/(2+M^{-2})$,
 a bound that is (often) poor, 
 since $M$ is usually much larger than $N$. 
 For instance, for the previous (binomial-type) examples 
 one finds the inaccurate bounds $R_3'\leq 5/11$ ($M=4$)
 and $R_3'\leq 85/171=0.497...$ \ ($M=16$), 
 respectively.
 
 Finally, fix $N\geq 2$ and define
 \[
 {\cal M}_N:=\Big\{{\bbb p}=(p_1,\ldots,p_N): 0\leq p_k<1 \ \ (k=1,\ldots,N), \ \
 p_1+\cdots+p_N=1\Big\}.
 \]
 It is fairly expected that, as ${\bbb p}$ varies in ${\cal M}_N$,
 the coefficient of correlation
 of any two order statistics from ${\cal P}_0({\bbb p})$
 will be maximized in the uniform case,
 ${\bbb p}=(1/N,\ldots,1/N)$, but no proof is available; we also conjecture
 that a similar result is true for the without-replacement case.

 \end{document}